\documentstyle[12pt]{article}

\newcommand{\intt}[1]{\rule{#1mm}{0mm}}
\def\ardl{ \mbox{\LARGE $\swarrow$}}
\def\ardr{ \mbox{\LARGE $\searrow$}}
\def\hv0{ \raisebox{2mm}{\mbox{\scriptsize $h \rightarrow 0$}}}
\def\tv0{ \raisebox{2mm}{\mbox{\scriptsize $t \rightarrow 0$}}}
\def\tetav0{ \raisebox{2mm}{\mbox{\scriptsize $\theta \rightarrow 0$}}}
\def\tauv0{ \raisebox{2mm}{\mbox{\scriptsize $\tau \rightarrow 0$}}}
\def\uat{ U(A_t)}
\def\upa{ U_p(A)}

\def\uhat{ U_h (A_t)}
\def\funhgt{ \mbox{Fun}_h(G_t)}
\def\funtgh{ \mbox{Fun}_t(G_h)}
\def\fungt{ \mbox{Fun}(G_t)}
\def\fungh{ \mbox{Fun}(G_h)}
\def\fungsh{ \mbox{Fun}(G^*_h)}
\def\funhab{ \mbox{Fun}_h ({\cal AB})}
\def\funab{ \mbox{Fun} ({\cal AB})}
\def\uhab{ U_h(\mbox{Ab})}
\def\uash{ U(A^*_h)}
\def\uast{ U(A^*_t)}

\begin{document}

\title{Classical limits,
quantum duality and Lie-Poisson structures}
\author{V.D. Lyakhovsky \thanks{Theoretical Department,
Institute of Physics,
St. Petersburg State University,
198904 St. Petersburg,
Russia,
e-mail: lyakhovsky@phim.niif.spb.su}
\\
  Institute of Theoretical Physics \\
  University of Wroc{\l}aw\\
  Pl. Maxa Borna 9\\
  50-204  Wroc{\l}aw \\
  Poland}
\date{}
\maketitle

\begin{abstract}
Quantum duality principle is applied to study classical limits
of quantum algebras and groups.  For a certain type of Hopf algebras
the explicit procedure to construct both classical limits is
presented. The canonical forms of quantized
Lie-bialgebras are proved to be two-parametric varieties with two
classical limits called dual. When considered from the point of view of
quantized symmetries such varieties can have boundaries that are
noncommutative and noncocommutative. In this case the quantum duality
and dual limits still exist while instead of Lie bialgebra one has a
pair of tangent vector fields. The properties of these
constructions called quantizations of Hopf pairs are studied and
illustrated on examples.
\end{abstract}

        \section{Introduction}
        Quantum duality principle \cite{Drin,STS} asserts
        that quantization of a Lie bialgebra $ (A,A^*) $
        gives rise to a dual pair of Hopf algebras
        $ (U_p (A),U_p (A^*)) $ or in dual terms --
        ($ \mbox{Fun}_p (G)$, $\mbox{Fun}_p (G^*) $).
        In the standard form
        quantum algebras and groups do not exhibit this
        duality explicitly. This is clearly seen when the
        classical limit is concerned. The quantum algebra
        $ U_p(A) $ is meant as a member of 1-dimensional
        family of Hopf algebras -- the deformation curve.
        Its classical limit $ U(A) $ is a fixed point of the
        orbit $ \mbox{Orb}(U(A)) $ where the deformation curve
        starts. Due to quantum duality formulated in terms
        of quantum formal series Hopf algebras \cite{Drin}
        quantum algebra can be interpreted as a quantum group
        \[
        U_p(A) \approx (\mbox{Fun}(G^*))_{p}
        \]
        for the universal covering group $G^* $ with Lie
        algebra $ A^* $. So there must be another classical
        limit, i.e. another deformation curve that starts
        at a fixed point of $ \mbox{Orb}(\mbox{Fun}(G^*)) $ and
        contains the Hopf algebra $ U_p(A) $.

        Thus the natural form of deformation quantization
        of Lie bialgebra $ (A,A^*) $ must be a 2-parametric
        family of Hopf algebras with two dual classical limits.
        Within certain assumptions this family forms an
        analytic variety $Q$ and the classical limits -- its boundary.
        The existence of a variety $Q$ with such properties is equivalent
        to attributing its member the quantum duality.
        The Lie bialgebra appears here in the form of two vector
        fields tangent to $Q$. From the point of view of Lie-Poisson
        structures and their symmetries it can be shown
        natural to consider the varieties of this type and their
        boundaries entirely placed in the domain of noncommutative
        and noncocommutative Hopf algebras. Preserving the main
        property of quantum duality we find that in this case some
        other characteristics are not conserved. In particular the
        lifted tangent fields may not form a Lie-bialgebra any more.
        But dual parameters and dual limits are still present there and
        can play important role in applications.

        The paper is organised as follows. In the subsection \ref{gen21}
        we describe how to construct the dual classical limits and under
        what conditions it can be done. The explicit example is considered
        in subsection \ref{ex22} where it is also shown how inverting
        the dual variety $Q^*$ one can change the role of dual
        parameters. In the subsection \ref{sub31}
        lifted varieties $Q_{\varepsilon}$ (that are called the quantized
        Hopf pairs) are studied. Their Lie-Poisson properties are
        considered in \ref{LPstr} and in \ref{ntexamp} the nontrivial
        example of such $Q_{\varepsilon}$ variety is presented. In
        Appendix the two-parametric form for the standard quantization
        of $sl(n,{\bf C})$ is given explicitly.

        \section{Dual classical limits}
        \subsection{General scheme} \label{gen21}
        Let us construct the second classical limit for a quantum Lie
        algebra $ U_p(A) $. Consider the variety $ {\cal H} $ of Hopf
        algebras with fixed number of generators. Its points
        $ H \in {\cal H} $ are parametrized by the corresponding
        structure constants. We must find  a (smooth)
        curve in $ {\cal H} $ containing $ U_p(A) $ and
        intersecting with the orbit $ \mbox{Orb}(\mbox{Fun}(G^*)) $.
        In the limit to be obtained the multiplication
        $ m $ in $ U_p(A) $ must become Abelian. For the
        universal enveloping algebra $ U(A) $ such a procedure
        is trivially described by a linear contraction.
        The corresponding transformation of basis $ \{ a_i \} $,
        \begin{equation}
        B(t) : a_i \rightarrow a_i/t ,
                                        \label{e1}
        \end{equation}
        leads to new structure constants $ C_{ij}^{'k} $
        \begin{equation}
        C_{ij}^{'k} (t) = t C_{ij}^{k}.
                                        \label{e2}
        \end{equation}
        The costructure  of $ U(A) $ being primitive is
        insensitive to this transformation. Algebras
        $ U(A_t) $ form a line in $ \mbox{Orb} U(A) $ with the
        limit point  $ U(A_0 \equiv \mbox{Abelian}) $. Let us apply
        operators $ B(t) $ to $ U_p(A) $. Such a sequence
        of equivalence transformations (for each value
        of $ p $ fixed) generates the smooth one-parametric
        curve
        \begin{equation}
        B(t) U_p(A) \equiv U_p (A_t)
                                        \label{e3}
        \end{equation}
        belonging to the orbit $ \mbox{Orb}(U_p(A)) $.
        Note that in general case Hopf algebras $ U_p(A) $
        are not equivalent for different $ p $. We thus
        obtain in $ {\cal H} $ the 2-parametric subset
        \begin{equation}
        \{ U_p (A_t) \}_{p>0,t>0} \equiv Q(A,A^*)
                                        \label{e4}
        \end{equation}
        formed by the dense set of smooth curves.
        Being the deformation quantization
        $ \{ U_p(A) \}_{p>0} $ is an analytic family.
        So $ Q(A,A^*) $ is a smooth variety with the global
        co-ordinates $p$ and $t$. Nevertheless these
        co-ordinates are not appropriate for our task
        because the limit
        \begin{equation}
        lim_{t\rightarrow 0} \, \, U_{p>0} (A_t)
                                        \label{e4a}
        \end{equation}
        does not exist. This is clearly seen from the
        properties of the coproduct  in $ U_p(A) $,
        \begin{equation}
        \Delta (a_i) = a_i \otimes {\bf 1} + {\bf 1} \otimes a_i +
        p (D_{i}^{kl} a_k \otimes a_l + S_{i}^{kl} a_k
        \otimes a_l + D_{i}^{ \{k,m\} \{l\} } a_k  a_m
        \otimes a_l + \cdots ).
                                        \label{e5}
        \end{equation}
        Here the first deforming function is divided
        into symmetric and antisymmetric parts and
        $ D_{i}^{kl} $ are the structure constants of $ A^* $.
        The transformation (3) leads to the expression
        divergent for $ t \rightarrow 0 $:
        \begin{equation}
        \Delta (a_{i}) = a_{i} \otimes {\bf 1} + {\bf 1} \otimes
        a_{i} + p/t (D_{i}^{kl} a_{k} \otimes a_{l}
         + \cdots ).                      \label{e6}
         \end{equation}
        Here the unwritten terms in (\ref{e6}) may contain the
        higher negative powers of $t$.

        It is easy to define a class of deformation quantizations
        for which one can overcome this difficulty. For each $ U_p(A) $
        consider the family of groups parametrized by $ p $ and defined
        on the space $ \mbox{Mor} (U_p(A),{\bf K}) $
        by the convolution multiplication
        \begin{equation}
        \phi_1 * \phi_2 = ( \cdot_{\bf K} )( \phi_1 \otimes
                          \phi_2 ) \Delta, \rule{1cm}{0cm}
        \phi_1, \phi_2 \in  \mbox{Mor} (U_p(A),{\bf K})
                                        \label{e6a}
        \end{equation}
        (${\bf K}$ here is the main field).
        Let ${\cal P}$ be the subspace
        of functionals dual to the space of algebra $A$.
        On this space (when certain conditions are fulfilled \cite{Lya})
        the convolution (\ref{e6a}) also generates a family of groups
        parametrised by $p$.
        This analytic family describes the
        contraction of the obtained group ${\cal P}$ to the additive
        Abelian vector group on the space dual to that of Lie algebra $A$.

        Suppose now that
        \begin{description}
        \item[(a)] the system of equations
\[  ( \cdot_{U_p} )(\mbox{id} \otimes S) \Delta =
( \cdot_{U_p} )(\mbox{id} \otimes S) \Delta = \eta \varepsilon
\]
on the basic elements $a_i$ fixes the antipode $S$ of $\upa$
(see \cite{Lya}) .
        \item[(b)] for $U_p(A)$
the contraction ${\cal P_{\rm p}} \longrightarrow {\cal AB}$
is equivalent to the trivial, that is induced by such a contraction of
its algebra $A^*$ where the structure constants in the lowest order
are proportional to $p$.
        \end{description}
        It is easy to verify that for such deformation quantizations
        in every monome of the coproduct (\ref{e5}) the power of $p$ is
        less than the total degree
        of the basic elements $a_i$ by one.

        If the Hopf algebra $U_p(A)$ belongs to the class described above
        the second classical limit can be obtained as follows.
        Let us change the co-ordinates:
        \begin{equation}
        (p,t)  \Rightarrow  (h,t),   \rule{1cm}{0cm}
        h  =  p/t.                        \label{e7}
        \end{equation}
        Then the coproduct (\ref{e6}) becomes well defined in the limit
        $t \rightarrow 0 $ (with $h$ fixed)
        because now its structure constants can bare only positive
        powers of $t$ (and in the two lowest orders can depend only on $h$).

        The multiplication
        structure constants in these new co-ordinates also have
        the finite limit values. To see this let us go to the
        dual picture.
        Consider the space of linear functionals on the
        space of Hopf algebra  $ U_p (A_t) $ and the canonical
        dual basis $ \{f^{I}\} $ such that the elements
        $ f^i $ dual to the basic elements of $A$,
        \[ <a_i, f^j >  = \delta_{i}^{j} , \]
        form the basis of $ A^* $.
        Construct the dual Hopf
        algebras for the elements of $ Q(A,A^*) $
        \begin{equation}
        \{ (U_p (A_t))^* \}_{p>0,t>0} =
        \{ \mbox{Fun}_p (G_t) \}_{p>0,t>0} \equiv
        Q^* (A,A^*)
                                        \label{e8}
        \end{equation}
        and also for the border line
        \begin{equation}
        \{ (U(A_t))^* \}_{t\geq0} = \{ \mbox{Fun}(G_t) \}_{t\geq0}.
                                        \label{e9}
        \end{equation}
        We shall consider these Hopf algebras as quantum
        formal series groups \cite{Wor}. The basis transformation
        $ B^*(t) $ dual to (\ref{e1})
        \begin{equation}
        B^* (t) : f^i \rightarrow t f^i
                                        \label{e10}
        \end{equation}
        leads to the following compositions:
        \begin{equation}
        [ f^i,f^k ] = p/t ( D_{l}^{ik} f^l + tE_{su}^{ik} f^s f^u
        + \cdots ),
                                        \label{e11}
        \end{equation}
        and
        \begin{equation}
        \Delta f^i = f^i \otimes {\bf 1} + {\bf 1} \otimes f^i +
        t (C_{kl}^{i} f^k \otimes f^l + T_{kl}^{i} f^k
        \otimes f^l + t C_{\{k,l\}\{u\}}^{i} f^k f^l
        \otimes f^u + \dots )
                                        \label{e12}
        \end{equation}
        For simplicity the renormalization factors for high
        power monomials are incorporated in the corresponding
        structure constants. These coproduct structure constants
        (in general) depend on $p$
        but they all have well defined limits when
        $ p \rightarrow 0 $. These limits describe
        the power series expansion of the
        multiplication law in $ G_t $ in terms of exponential
        co-ordinate functions. So the corresponding Taylor
        series can be written for them in the neighbourhood
        of $ ((\cdot))p = 0 $. Substituting $ p=ht $ in (\ref{e12})
        and going to the limit $ t \rightarrow 0 $ we see
        that all coefficients are finite. (Note that in the
        framework of formal series Hopf algebras this
        conclusion is true for all $ p \geq 0 $.)

        We have shown that for algebras $ U_h (A_t) \in
        Q (A,A^*) $ of the described class the limit $ U_h (A_0) \equiv
        lim_{t\rightarrow 0} \, \, U_h (A_t) $  exists. According
        to the quantum duality
        $ U_h (A_t) \approx \mbox{Fun}_t (G_{h}^{*}) $,
        so
        \begin{equation}
        lim_{t \rightarrow 0} \, \,  U_h (A_t) \equiv U_h (A_0)
        \approx \mbox{Fun} (G_{h}^{*})
                                        \label{e13}
        \end{equation}
        Thus every such deformation quantization can be written
        in the form $ U_h (A_t) $ (respectively
        $ \mbox{Fun}_t (G_{h}^{*}) ) $ that reveals two canonical
        dual classical limits:
\begin{equation}
\begin{array}{cc}
\intt{5} \uhat \intt{5} & \\
\hv0 \ardl \intt{7} \Updownarrow \intt{7} \ardr \tv0 & \\
\intt{3} \uat \intt{3} \intt{10} \funhgt \intt{10} \fungsh & \approx \uhab \\
\Updownarrow \intt{8} \ardl \intt{20} \ardr \intt{8} \Updownarrow &  \\
\intt{2} \fungt \intt{38} \funhab & \approx \uash
\end{array}                 \label{e14}
\end{equation}

        All the reasoning is invariant with respect to
        interchange $ A \rightleftharpoons A^* $.
        So the set $ \{ \mbox{Fun} (G_{h}^{*}) \} $ can be also
        considered as a straight line in the orbit
        $ \mbox{Orb} (\mbox{Fun}(G^*)) $ -- the trivial contraction of
        the group $ G^* $ into the abelian additive vector
        group $ {\cal AB} $. The lines $ \{ U (A_t) \} $ and
        $ \{ \mbox{Fun} (G_{h}^{*}) \} $ intersect in the point
                $ U(\mbox{Ab}) \approx \mbox{Fun} ({\cal AB}) $.

The conclusion can be formulated as follows.
\newtheorem{propo}{Proposition}
\begin{propo}
If $\upa$ is a deformation quantization of a Lie bialgebra $(A,A^*)$
with the properties (a)-(b) then in ${\cal H}$ there exists an analytic
submanifold $Q(A,A^*)$ (respectively $Q^*(A,A^*)$).It can be globally
parametrised by co-ordinates $(h,t)$. The trivial contraction lines
$\uat$ and $\fungsh$ (respectively $\uash$ and $\fungt$) together with
their intersection point $U(\mbox{Ab}) = \funab $ form a boundary of
$Q(A,A^*)$ (respectively $Q^*(A,A^*)$) and supply its elements with
the dual classical limits:
         \begin{equation}
         \left.
         \begin{array}{l}
         lim_{h \rightarrow 0} H_{h,t} = U (A_t)  \\
         lim_{t \rightarrow 0} H_{h,t} = \mbox{Fun} (G_{h}^*)
         \end{array}
         \right\} \mbox{for} \, \, Q            \label{e15}
         \end{equation}
         and
         \begin{equation}
         \left.
         \begin{array}{l}
         lim_{h \rightarrow 0} H_{h,t}^{*} = \mbox{Fun} (G_t)  \\
         lim_{t \rightarrow 0} H_{h,t}^{*} = U (A_{h}^{*})
         \end{array}
         \right\} \mbox{for} \, \, Q^* .        \label{e16}
         \end{equation}
\end{propo}
        Suppose now that
        \begin{enumerate}
        \item $ t $ and $ h $ parametrize
        the intersecting trivial contraction lines
        $ U (A_t) $ and $ \mbox{Fun} (G_{h}') $ in the orbits
        $ \mbox{Orb} (U(A)) $ and $ \mbox{Orb} (\mbox{Fun}(G')) $
        (respectively
        $ \mbox{Orb} (\mbox{Fun}(G)), \mbox{Orb} (U(A')) $ )
        for a certain pair of inequivalent algebras $ (A,A^{'}) $
        of equal dimention and corresponding universal covering
        Lie groups $G$ and $G'$.
        The intersection
        point coincides with the contraction limit.
        \item In the variety
        $ {\cal H} $ of Hopf algebras with the generators
        $ \{ a_i, {\bf 1} \} $ ($ \{ a_i \} $ is the basic of
        $ A $) there exists the analytical 2-dimensional
        subvariety $ Q (A,A') $ (respectively
        $ Q^* (A,A') $) of Hopf algebras $ H $ and the
        disunqued union $ U(A_t) \cup \mbox{Fun} (G_{h}') $
        is the boundary of $ Q(A,A') $.
        \end{enumerate}
It is easy to check the validity of the following statement
\begin{propo}
If for algebras $ A $ and $ A' $ the conditions (1)-(2) are fulfilled
-- they are dual,
Hopf algebras $ H \in Q(A,A') $ are the deformation
         quantizations of the  Lie bialgebra
         $ (A,A' \approx A^* ) $ and the contraction curves
$\uat$ and $\fungsh$ supply the dual classical limits to the points of
$Q(A,A^*)$.
\end{propo}

         It is possible to invert the dual list $ Q^* $
         with respect to $ Q $ so that the dual limits
         will refer to dual algebras (respectively groups).
         This is due to the fact  that in the
         $ (h,t) $ - co-ordinates the Hopf algebras
         $ H_{h,t}^{*} $ and $ H_{t,h}^{*} $ are equivalent,
         they are characterised by the same parameter
         $ p = ht $ (see (\ref{e7})) and are connected
         by the transformation $ B^* (h/t) $.
        Thus one can always introduce the transformed nondegenerate
        bilinear form $ < , >_{h/t} $ such that the parameters
        on the list $ Q^* $ will be interchanged with respect to
        $ Q $,
        \begin{equation}
        H_{h,t} \stackrel{ <>h/t }{ \Longleftrightarrow } H_{t,h}^* .
                                        \label{e17}
        \end{equation}
        Now the correlation between the canonical classical limits
        for $ Q $ and $ Q^* $ will differ from
        that described by the diagram (\ref{e14}):
\begin{equation}
\begin{array}{c}
\intt{5} \uhat \intt{5}  \\
\hv0 \ardl \intt{2}
\Updownarrow <> \mbox{\scriptsize $h/t$}  \intt{2} \ardr \tv0  \\
\intt{4} \uat \intt{2}
\intt{3} \Updownarrow \intt{3} \funtgh \intt{3} \Updownarrow \intt{3}
\fungsh  \\
 \hv0 \ardl \intt{20} \ardr \tv0 \\
\intt{3} \uast \intt{37} \fungh
\end{array}                      \label{e18}
\end{equation}
        Note that here in contrast to (\ref{e14}) the duality holds only
        for points of $ Q $ and $ Q^* $.

\subsection{Example: $ U_h(sl_t(2)) \approx
        \mbox{Fun}_t (\widetilde{E}_h(2))$ }    \label{ex22}
        Consider the standard quantum algebra $ U_p(sl(2,{\bf C})) $:
        \begin{equation}
        \begin{array}{lll}
        \Delta L = L \otimes {\bf 1} + {\bf 1} \otimes L, & \rule{5mm}{0mm} &
        [L,M] = M,                                        \\
        \Delta N = N \otimes {\bf 1} + e^{-2pL} \otimes N, &  &
        [L,N] = -N,                                       \\
        \Delta M = M \otimes e^{2pL} + {\bf 1} \otimes M, &  &
        [M,N] = 2 \frac{\sinh 2pL}{1 - e^{-2p}}.
        \end{array}                     \label{e19}
        \end{equation}
        Applying the transformation $ B(t) $ and introducing the new
        co-ordinates $ (h=p/t, t) $ (see (\ref{e7})) we get the canonical
        form for Hopf algebras of the variety $ Q(sl(2), \widetilde{e(2)})$:
        \begin{equation}
        \begin{array}{lll}
        \Delta L = L \otimes {\bf 1} + {\bf 1} \otimes L, & \rule{5mm}{0mm} &
        [L,M] = tM,                                        \\
        \Delta N = N \otimes {\bf 1} + e^{-2hL} \otimes N, & &
        [L,N] = -tN,                                       \\
        \Delta M = M \otimes e^{2hL} + {\bf 1} \otimes M, & &
        [M,N] = 2t^2 \frac{\sinh 2hL}{1 - e^{-2ht}}.
        \end{array}                     \label{e20}
        \end{equation}
        The dual classical limits are
        \begin{equation}
        \begin{array}{c}
        \lim_{h \rightarrow 0} U_h(sl_t(2)) = U(sl_t(2)), \\
        \lim_{t \rightarrow 0} U_h(sl_t(2)) = \mbox{Fun} (\widetilde{E_h(2)}).
        \end{array}                     \label{e21}
        \end{equation}
        Here $ \widetilde{E(2)} $ is the group of flat motions with
        the quasiorthogonal rotation.

        The canonical dualization produces the variety $ Q^*(sl(2),
        \widetilde{e(2)}) $ whose points can be interpreted as
        $ \mbox{Fun}_h(SL_t(2)) $ and the parametrization is as follows:
        \begin{equation}
        \begin{array}{rcl}
       \rule{1mm}{0mm} [\lambda,\mu] & = & -2h \mu, \\
       \rule{1mm}{0mm} [\lambda,\nu] & = & -2h \nu, \\
       \rule{1mm}{0mm} [\mu, \nu]    & = & 0,       \\
        \Delta \mu    & = & {\bf 1} \otimes \mu + ( \mu \otimes e^{-t\lambda})
                ({\bf 1} \otimes {\bf 1} + t^2 \mu \otimes \nu )^{-1} \\
        \Delta \nu    & = & \nu \otimes {\bf 1} + ( e^{-t\lambda} \otimes \nu)
                ({\bf 1} \otimes {\bf 1} + t^2 \mu \otimes \nu )^{-1} \\
        \Delta \lambda &= & \lambda \otimes {\bf 1} + {\bf 1} \otimes \lambda
        -4ht\sum_{n=1}^{\infty} t^{2n-1}\frac{(- \mu \otimes \nu )^n }{ 1-
        e^{-2nht}}.
        \end{array}                     \label{e22}
        \end{equation}
        We have used here the realization of $ \mbox{Fun}_h (SL(2)) $
        obtained in \cite{CELL}. On the dual list $ Q^* $ the classical
        limits are
        \begin{equation}
        \begin{array}{c}
        \lim_{h \rightarrow 0} \mbox{Fun}_h(SL_t(2)) = \mbox{Fun}(SL_t(2)), \\
        \lim_{t \rightarrow 0} \mbox{Fun}_h(SL_t(2)) = U (\widetilde{e_h(2)}).
        \end{array}                     \label{e23}
        \end{equation}
        The first of these two limits is treated here as
        the formal series group whose generic elements are
        the basic exponential co-ordinate functions on
        $ SL_t(2) $ in the neighbourhood of unit.

        Using the equivalence transformation $ B^*(h/t) $ on $ Q^* $
        the parametric dualization can be defined:
        \begin{equation}
        <L,\lambda>_{h/t} =<M,\mu>_{h/t} =<N,\nu>_{h/t} = h/t.
                                        \label{e24}
        \end{equation}
        It inverts the dual list $ Q^* $ with respect to $ Q $
        and shows explicitly that quantum algebras of $A$ and $A^*$
        as well as quantum groups of $G$ and $G^*$ form the dual pairs of
        Hopf algebras
        (see (\ref{e18})):
\begin{equation}
\begin{array}{c}
\intt{5}  U_h (sl_t(2)) \intt{5}    \\
\hv0 \ardl \intt{8}
\Updownarrow <> \mbox{\scriptsize $h/t$}
\intt{8} \ardr \tv0 \\
\intt{4} U(sl_t(2)) \intt{1} \Updownarrow
\intt{4} \mbox{Fun}_t(SL_h(2)) \intt{4}  \Updownarrow
\mbox{Fun} (\widetilde{E_h(2)})   \\
\hv0 \ardl \intt{30} \ardr \tv0   \\
\intt{6} U(\widetilde{e_t(2)}) \intt{44} \mbox{Fun} (SL_h(2))
\end{array}                      \label{e25}
\end{equation}

        In the Appendix the 2-parametric variety $ Q $ for quantum algebras
        $ U_h(sl(n,{\bf C})) $ is explicitly described.

        \section{Quantization of Hopf pairs}
        \subsection{ } \label{sub31}
        The geometric description of the quantum duality properties of the
        deformation quantization given by the criteria (1)-(2) leads to the
        natural generalisation of this notion.

        Let us call the deformation quantization of the Hopf pair
        $ (H_{(\theta, 0)},H_{(0,\tau)}) $ the following construction in
        the variety $ {\cal H} $ of Hopf algebras with the generators
        $ \{ a_i, {\bf 1} \} $ : \\
        (1') $ \theta $ and $ \tau $ parametrize the intersecting
        contraction curves $ H_{(\theta,0)} $ and $ H_{(0,\tau)} $ such that
        in the first order the coproducts of generators remain undeformed
        in $ H_{(\theta,0)}$ and the products of generators
        -- in $ H_{(0, \tau)}$;
        the intersection point $ H_{(0,0)} $
        coincides with the contraction limit for both
        $ H_{(\theta, 0)} $ and $ H_{(0,\tau)} $. \\
        (2') In the variety $ {\cal H} $ there exists the analitic
2-dimensional
        manifold $ Q(H_{(\theta, 0)},H_{(0,\tau)}) $ (respectively
        $ Q^* (H_{(\theta, 0)},H_{(0,\tau)}) $ of (noncommutative and
        noncocommutative) Hopf algebras $ H_{(\theta,\tau)} $ . The
        disjoint union
        $ H_{(\theta, 0)} \bigcup H_{(0,\tau)} \bigcup H_{(0,0)} $
        is the boundary of $ Q(H_{(\theta, 0)},H_{(0,\tau)}) $.
        (Respectively
        $ H^*_{(\theta, 0)} \bigcup H^*_{(0,\tau)} \bigcup H_{(0,0)} $ --
        the boundary of $ Q^* $ .) \\
        Then it is easy to check that on the curves
        $ H_{(\theta,0)} $ and $ H_{(0,\tau)} $
        the smooth vector fields of the first deforming functions
        $ V(\theta) $ and $ W(\tau) $ exist such that the limit vectors
        \begin{equation}
        \lim_{\theta \rightarrow 0} V(\theta) \equiv V(0), \label{e26}
        \end{equation}
        \begin{equation}
        \lim_{\tau \rightarrow 0} W(\tau) \equiv W(0), \label{e27}
        \end{equation}
        are the first deformation functions for $H_{(0,0)}$
        in the direction of $ H_{(0,\tau)} $ and $ H_{(\theta,0)} $
        respectively.

        To connect this purely geometric construction with the usual
        deformation quantization picture it is natural to formulate
        the additional condition. \\
        (3') In $ {\cal H} $ there exists a smooth 1-parametric family
        $ Q_{\varepsilon} $ of varieties
        $ Q(H_{(\theta, 0)},H_{(0,\tau)}) $ such that
        its limit
        \[
        \lim_{\varepsilon  \rightarrow 0} Q_{\varepsilon} \equiv
                 Q_0(H_{0(\theta, 0)},H_{0(0,\tau)})
        \]
        satisfies the conditions (1)-(2).

        The main difference with the canonical case is that the
        intersection point algebra $ H_{\varepsilon} (0,0) $ may be
        noncommutative and noncocommutative for $ \varepsilon \neq 0 $.
        From this point of view
        the canonical dualization is the dualization with respect to
        the Abelian and coAbelian Hopf algebra. While in the
        quantization of a Hopf pair no such restrictions on
        $ H_{\varepsilon} (0,0) $ are imposed.
        The parameters
        $\theta$ and $\tau$ are dual with respect to the common limit
        $ H_{0,0} $ of the boundary curves $ H_{(\theta,0)} $
        and $ H_{(0,\tau)} $.
        The pair of vector fields
        $ (V_{\varepsilon}(\theta),W_{\varepsilon} (\tau)) $
        plays here the role of a Lie-bialgebra.
        In the limit $ \varepsilon \rightarrow 0 $ the pair
         $ (V_0(0), W_0(0)) $
        becomes a Lie-bialgebra.

        The existence of the deformation quantization of a Hopf pair
        is tightly connected with the contraction properties of Hopf
        algebras. In \cite{BAL} it was demonstrated how the deformation
        parameters  of the quantum group can be dualized with the
        quantization parameters of the corresponding quantum algebra.
        The algebra $ U_{\phi} {\cal G}_{0,k_2} $ was constructed to
        illustrate this effect:
        \begin{equation}
        \begin{array}{l}
       \rule{1mm}{0mm} [P_3,P_1] =
        \frac{e^{i \phi P_2} - e^{-i \phi P_2}}{2i \phi}, \\
       \rule{1mm}{0mm} [P_3,P_2] = -k_2 P_1, \\
       \rule{1mm}{0mm} [P_1,P_2] = 0.
        \end{array}                     \label{e28}
        \end{equation}
        \begin{equation}
        \begin{array}{rcl}
        \Delta(P_{1,3})& = &e^{-i \phi P_2 /2} \otimes P_{1,3} +
                        P_{1,3} \otimes  e^{i \phi P_2 /2}, \\
        \Delta(P_2)& = & P_2 \otimes {\bf 1} + {\bf 1} \otimes P_2.
        \end{array}                             \label{e29}
        \end{equation}
        The authors emphasise that in the dual Hopf algebra
        $ (U_{\phi} {\cal G}_{0,k_2})^* $ the contraction parameter
        $ k_2 $ measures the deformation of coproduct and thus obtain the
        features of the quantization parameter.

        From our point of view it must be also stressed that in the limit
        \[
        \lim_{k_2 \rightarrow 0} U_{\phi} {\cal G}_{(0,k_2)} =
        U_{\phi} {\cal G}_{(0,0)}
        \]
        we do not get the classical algebra of functions on the dual
        group . The co-ordinate functions remain noncommutative. So the
        parameters $ k_2 $ and $ \phi $ are not canonically dual.

        The situation becomes clear in terms of
        quantization of a Hopf pair. The contraction curves
        $ \mbox{Fun}_{\phi} (G^*_{\phi} ) \equiv H_{(0,\phi)} $ and
        $ U {\cal G}_{(0,k_2)} \equiv H_{(k_2,0)} $
        have the common limit
        $ U {\cal G}_{(0,0)} \equiv H_{(0,0)} $. Here $H_{(0,0)}$ is
        the Hiesenberg algebra. The union
        $ \mbox{Fun}_{\phi} (G^*_{\phi} ) \bigcup U {\cal G}_{(0,k_2)} $
        form the boundary of the 2-dimensional variety
        $ Q(H_{(k_2,0)},H_{(0,\phi)}) $ . The vector fields $ V $ and $W$
        are trivial in the sense that their projections
$ V_{\downarrow {\cal G} \wedge {\cal G} \rightarrow {\cal G}} $
        and
$ W_{\downarrow {\cal G} \rightarrow {\cal G} \wedge {\cal G}} $
do not depend on
        $ k_2 $ and $ \phi $ respectively. The first deforming function
        \[
        [P_3,P_2]_{(1)} = -P_1
        \]
        with the only nonzero co-ordinate $ V^1_{[3,2]} = -1 $ generates
        \[
        V(\phi) = V(0)
        \]
        and the first deforming function
        \[
        \Delta_{(1)} P_{1,3} = -i/2(P_2 \otimes P_{1,3} - P_{1,3} \otimes
                                P_2 )
        \]
        with nontrivial co-ordinates $ W_1^{[2,1]} = -i/2, \, \,
        W_3^{[2,3]} = -i/2 $ generates the field
        \[
        W(k_2) = W(0).
        \]
        Notice that here the pair
$ (V_{\downarrow {\cal G} \wedge {\cal G} \rightarrow {\cal G}},
 W_{\downarrow {\cal G} \rightarrow {\cal G} \wedge {\cal G}}) $
is just a Lie-bialgebra.
        This fact is tightly connected with the possibility to fulfil
        the condition (3'). The family $ Q_{\varepsilon} $ is
        obtained by a simple transformation of generators
        \[
        P_{1,3} = p_{1,3}/ \varepsilon, \, \, \, P_2 = p_2
        \]
        that doesn't touch the tangent vectors $(V,W)$.
        In terms of new generators $ p_i$ the left-hand side of the first
        commutator in (\ref{e28}) acquires the multiplier $ \varepsilon^2 $
        and in the limit $ \varepsilon \rightarrow 0 $ the 2-dimensional
        subvariety $ Q_0 $ will have all the properties characteristic
        to the canonical deformation quantization scheme (see
        the conditions (1)-(3)):
        \begin{equation}
        \begin{array}{llcl}
       \rule{1mm}{0mm} [p_3,p_1] = 0, &
        \Delta(p_{1,3})& = &e^{-i \phi p_2 /2} \otimes p_{1,3} +
                        p_{1,3} \otimes  e^{i \phi p_2 /2}, \\
       \rule{1mm}{0mm} [p_3,p_2] = -k_2p_1, &
        \Delta(p_2)& = & p_2 \otimes {\bf 1} + {\bf 1} \otimes p_2, \\
       \rule{1mm}{0mm} [p_1,p_2] = 0. & & &
        \end{array}                             \label{e30}
        \end{equation}
        We come to the conclusion that the variety described by
        (\ref{e28}) and (\ref{e29}) is the deformation
        of a Hopf pair $ (H_{(k_2,0)},H_{(0,\phi)}) $. The parameters $ \phi $
        and $k_2 $ are dual with respect to $ H_{\varepsilon (0,0)} $
        -- the universal enveloping algebra of the Hiesenberg algebra.
        In the limit $ \varepsilon \rightarrow 0 $ they become
        canonically dual.

        It is obviously clear that considering $ U_{\phi} {\cal G}_{(0,k_2)} $
        with fixed $ k_2 $ just as the canonical deformation quantization
        of $ U {\cal G}_{(0,k_2)} $, one can easily find the parameter $ \psi $
        canonically dual to $ \phi $. The standard procedure described in
        section 1 leads to the following variety $ Q( {\cal G}_{(0,k_2)},
        \widetilde{e(2)}) $,
        \begin{eqnarray*}
        \rule{1mm}{0mm} [P_3,P_1] & = & \frac{\psi}{2i \phi}
                                     (e^{i \phi P_2} - e^{-i \phi P_2}), \\
        \rule{1mm}{0mm} [P_3,P_2] & = & - \psi k_2 P_1, \\
        \rule{1mm}{0mm} [P_1,P_2] & = & 0,  \\
        \Delta P_{1,3} & = & e^{-i/2 \phi P_2} \otimes P_{1,3} +
                           P_{1,3} \otimes e^{i/2 \phi P_2}, \\
        \Delta P_2 & = & P_2 \otimes {\bf 1} + {\bf 1} \otimes P_2.
        \end{eqnarray*}
        In the considered example of the quantization of a Hopf pair
        one of the contraction curves
        -- $ H_{\varepsilon (k_2,0)} $ -- remains canonical
        (with a primitive coproduct) even for $ \varepsilon \neq 0 $.
        The reason is that the contraction described by the parameter
        $ k_2 $ does not touch the coproducts of basic coordinate functions.
        In terms of the quantization of a Hopf pair this means
        that only one of the contractions, namely the $ \mbox{Fun}
        ( G^*_{\phi} ) $ is "lifted" by the $\varepsilon$-deformation
        nontrivially. In the subsection \ref{ntexamp} we demonstrate
        an example
        where both contraction lines are nontrivially lifted.

        \subsection{Lie-Poisson structures in a case of quantized pair}
        \label{LPstr}
        In the deformation quantization of a Lie bialgebra the
        quantum algebra $ U_h(A_t) $ refers to the initial
        Lie-Poisson structure as to the Poisson-Hopf algebra
        $ (U(A_t), A^*) $, where the Poisson comultiplication is
        described by a Lie structure of $ A^* $. In construction
        of a Hopf algebra $ U_h(A_t) $ the comultiplication is
        given a preference -- it is quantized first. The deformations
        of commutation relations play an auxiliary role.
        It starts when the first order deformation with respect to
        parameter $h$ is already performed in $ \Delta_{\downarrow A} $'s
        by $ A^* $.

        Considering $ U_h(A_t) $ as a quantum group
        $ \mbox{Fun}_t(G^*_h) $ one obtains the dual scenario.
        Now one has a Lie-Poisson group $ (\mbox{Fun}(G^*_h),A) $.
        The multiplication is quantized first and the deformation
        of the costructure is auxiliary. For example, when the group
        $ G^*_h $ is Abelian the noncommutative co-ordinates can
        always be introduced without any deformation of $ G^*_h $.

        Having the applications in mind one can consider such a
        Lie-Poisson structure $ A^* $ on $ \mbox{Fun}(G_t) $
        that do not form a Poisson-Hopf algebra with $ U(A) $,
        but nevertheless the simultaneous Hopf deformations of both
        multiplication and comultiplication exist. In this case
        the first order deforming function of $ U(A_t) $ will
        have the nontrivial multiplication constituents, that
        naturally may not form a Lie algebra themselves. They
        depend on $t$ and tend to zero when $ t \rightarrow 0 $
        if $t$ is still dual to $h$ with respect to $ U(\mbox{Ab}) $.
        This means that the desired Poisson properties can be
        imposed on the group $G$ only together with the symmetry
        deformation.
        In other words the Poisson structure
        on $ \mbox{Fun} (G_t) $ defined by $ A'$ does not perform
        the group $ G_t $ into a Lie-Poisson group but does it for a
        certain infinitesimal deformation $ (G_t)_h $.
        From the Hopf pair quantization point of view this case
        is based on such a deformation quantization where the field
        $ V(t) $ ( or both $ V(t) $ and $ W(h) $ ) has not only
        comultiplicative but also multiplicative nontrivial parts. The Lie
        bialgebra is reobtained in the limit $ t,h \rightarrow 0 $
        (see (\ref{e26}), (\ref{e27})).

        It must be noticed that the canonical scheme is stable
        (in the sense of generalisation described above)
        for semisimple algebras $A$ (or $A^*$).
        For them a first deforming function can be always set to
        zero by a similarity transformation of $A$. But for
        physical applications the nonsemisimple algebras play an
        important role and for them a deformation with a nontrivial
        first order may exist.

        The other possibility is to consider deformations of
        quantum algebras induced by deformations of classical ones
        \cite{BAL}. The example considered in the \ref{sub31}
        illustrates this case. There the Hopf algebra $ H_{(0,0)} $
        is nontrivial but classical Lie universal enveloping.
        In general the pair of deforming functions $ (V(0),W(0)) $
        must not form a Lie bialgebra in such a case.
        The Lie bialgebra structure must be obtained as a limit
        for $ \varepsilon \rightarrow 0 $ if the condition (3')
        is fulfilled.
        This is the case when the Poisson structure that form a
        Poisson-Hopf algebra with $A$ and gives rise to a global
        deformation is applicable also for a classical deformation
        $A_t$ with the analogous results.

        The general case described by the conditions (1')-(3')
        can be considered as a combination of these two possibilities.
        Here instead of Poisson structures on Lie group $G$ a first
        deforming function of some initial Poisson structure on $G$ appears.
        And instead of constructing a Lie structure for the initial
        multiplication described by $A$ and supposed to become Poisson
        we see a deformation problem for some already existing
        Lie structure, compatible with $A$. The details of this
        situation are explicitly demonstrated in the following subsection.

        \subsection{Example of nontrivial lifting}  \label{ntexamp}

        To obtain a nontrivial lifting for both multiplications and
        comultiplications one must chose quantum algebras that have
        nontrivial deformations. A favourable situation may
        be found in case of Hopf algebras obtained as Drinfeld doubles
        \cite{Drin}.

        Let us construct the double $D$ for the Hopf algebras
        $ U_h (sl_t(2,{\bf C})) $ and $ \mbox{Fun}_h(\widetilde{E_t(2)}) $
        (see (\ref{e20},\ref{e22})). To make the compositions more
        transparent The modified notations for the basic elements,
        \[
        \begin{array}{c}
        \{ L,M\intt{2},N\intt{2},\intt{1}\lambda,
        \intt{1}\mu\intt{2},\nu\intt{3} \} \\
        \Downarrow \\
        \{ L,X_+,X_-,H,Y_+,Y_- \}
        \end{array}
        \]
        make the compositions of the double more transparent:
        \begin{equation}
        \begin{array}{rcl}
        \rule{1mm}{0mm} [L,X_{\pm}] & = & \pm h X_{\pm}, \\
        \rule{1mm}{0mm} [X_+,X_-] & = & h^2 \frac{e^{tL} - e^{-tL} }
        {1- e^{-th}} ,
        \end{array}                     \label{e31}
        \end{equation}
        \begin{equation}
        \begin{array}{rcl}
        \rule{1mm}{0mm} [H,Y_{\pm}] & = & -t Y_{\pm}, \\
        \rule{1mm}{0mm} [Y_+,Y_-] & = & 0 ,
        \end{array}                     \label{e32}
        \end{equation}
        \begin{equation}
        \begin{array}{rcl}
        \rule{1mm}{0mm} [H,L] & = & 0, \\
        \rule{1mm}{0mm} [H,X_{\pm}] & = & tX_{\pm} \pm
        \frac{2th^2 }
        {1- e^{-th}} Y_{\mp} , \\
        \rule{1mm}{0mm} [Y_{\pm},L] & = & \pm h Y_{\pm}, \\
        \rule{1mm}{0mm} [Y_{\pm},X_{\pm}] & = & ( \mp e^{\pm tL}
        \pm e^{-hH} ) \\
        \rule{1mm}{0mm} [Y_{\pm},X_{\mp}] & = & \pm h^2 Y^2_{\pm}
        +( e^{-th} - 1 ) X_{\mp} Y_{\pm},
        \end{array}                     \label{e33}
        \end{equation}
        \begin{equation}
        \begin{array}{rcl}
        \Delta L & = & L \otimes {\bf 1} + {\bf 1} \otimes L, \\
        \Delta X_+ & = & X_+ \otimes e^{tL}
        + {\bf 1} \otimes X_+, \\
        \Delta X_- & = & X_- \otimes {\bf 1}
        + e^{-tL} \otimes X_-,
        \end{array}                     \label{e34}
        \end{equation}
        \begin{equation}
        \begin{array}{rcl}
        \Delta H & = & H \otimes {\bf 1} + {\bf 1} \otimes H
        -2t \sum_{n=1}^{ } \frac{(-h^2 (Y_- \otimes Y_{+}))^n}{1 - e^{-nth}},
\\
        \Delta Y_+ & = & Y_+ \otimes {\bf 1}
        + (e^{-hH} \otimes Y_+)({\bf 1} \otimes {\bf 1}
        + h^2 Y_{-} \otimes Y_{+})^{-1}, \\
        \Delta Y_- & = & {\bf 1} \otimes Y_-
        + (Y_{-} \otimes e^{-hH})({\bf 1} \otimes {\bf 1}
        + h^2 Y_{-} \otimes Y_{+})^{-1}.
        \end{array}                     \label{e35}
        \end{equation}
        The costructure of this quantum algebra describes the direct
        product $ \widetilde{E_t(2)} \times SL_h(2,{\bf C}) $.

        Consider a pair of new parameters $ ( \tau, \theta ) $ that describe
        the contractions
        \begin{equation}
        \begin{array}{rcl}
        SL_h(2,{\bf C}) &
        \stackrel{\mbox{\scriptsize
        contract}}{\tau \longrightarrow 0} & E_h(2,{\bf C}) \\
        sl_h(2,{\bf C}) &
        \stackrel{\mbox{\scriptsize
        contract}}{\theta \longrightarrow 0} & e_h(2,{\bf C})
        \end{array}                     \label{e36}
        \end{equation}
        lifted to the variety ${\cal H}$ of quantum algebras with six
        generators. After the necessary re\-para\-met\-ri\-za\-tion of the type
        described in section 1 we get the two-dimensional family
        $ D_{\tau ,\theta} $ (on this stage of construction the
        parameters $t$ and $h$ may be fixed).
        The elements of this variety differ from the initial Hopf algebra
        $D$ in the following compositions:
        \begin{equation}
        \begin{array}{rcl}
        \rule{1mm}{0mm} [X_+,X_-] & = & \tau \theta^2 h^2
        \frac{e^{tL} - e^{-tL} }{1- e^{-th}} ,    \\
        \rule{1mm}{0mm} [H,X_{\pm}] & = & tX_{\pm} \pm  \tau \theta
        \frac{2th^2 }{1- e^{-th}} Y_{\mp} , \\
        \rule{1mm}{0mm} [Y_{\pm},X_{\pm}] & = & \theta (
        \mp e^{\pm tL} \pm e^{-hH} ) \\
        \rule{1mm}{0mm} [Y_{\pm},X_{\mp}] & = &
        \pm \tau \theta  h^2 Y^2_{\pm}
        +( e^{-th} - 1 ) X_{\mp} Y_{\pm},
        \end{array}                     \label{e37}
        \end{equation}
        \begin{equation}
        \begin{array}{rcl}
        \Delta H & = & H \otimes {\bf 1} + {\bf 1} \otimes H
        -2t \sum_{n=1}^{ }
        \frac{(-\tau h^2 (Y_- \otimes Y_{+}))^n}{1 - e^{-nth}}, \\
        \Delta Y_+ & = & Y_+ \otimes {\bf 1}
        + (e^{-hH} \otimes Y_+)({\bf 1} \otimes {\bf 1}
        + \tau h^2 Y_{-} \otimes Y_{+})^{-1}, \\
        \Delta Y_- & = & {\bf 1} \otimes Y_-
        + (Y_{-} \otimes e^{-hH})({\bf 1} \otimes {\bf 1}
        + \tau h^2 Y_{-} \otimes Y_{+})^{-1}, \\
        \end{array}                     \label{e38}
        \end{equation}
        The smooth subvariety $ D_{\theta , \tau} $ contains two nontrivial
        contraction curves $ D_{\theta , 0} $ and $ D_{0, \tau} $
        having the noncommutative and noncocommutative Hopf algebra
        $ D_{0,0} $ as the common limit:
        \begin{equation}
        \begin{array}{rcl}
        \rule{1mm}{0mm} [L,X_{\pm}] & = & \pm h X_{\pm}, \\
        \rule{1mm}{0mm} [H,Y_{\pm}] & = & -t Y_{\pm}, \\
        \rule{1mm}{0mm} [H,X_{\pm}] & = & tX_{\pm}, \\
        \rule{1mm}{0mm} [Y_{\pm},X_{\mp}] & = &
        ( e^{-th} - 1 ) X_{\mp} Y_{\pm},
        \end{array}                     \label{e39}
        \end{equation}
        (All other pairs of generators commute.)
        \begin{equation}
        \begin{array}{rcl}
        \Delta L & = & L \otimes {\bf 1} + {\bf 1} \otimes L, \\
        \Delta X_+ & = & X_+ \otimes e^{tL}
        + {\bf 1} \otimes X_+, \\
        \Delta X_- & = & X_- \otimes {\bf 1}
        + e^{-tL} \otimes X_-,
        \end{array}                     \label{e40}
        \end{equation}
        \begin{equation}
        \begin{array}{rcl}
        \Delta H & = & H \otimes {\bf 1} + {\bf 1} \otimes H, \\
        \Delta Y_+ & = & Y_+ \otimes {\bf 1}
        + e^{-hH} \otimes Y_+, \\
        \Delta Y_- & = & {\bf 1} \otimes Y_-
        + Y_{-} \otimes e^{-hH}.
        \end{array}                     \label{e41}
        \end{equation}
        The picture of the smooth subvariety $ D_{\theta , \tau} $
        and its boundary looks like follows:
\begin{equation}
\begin{array}{c}
\ardr \intt{23} D_{\theta' ,\tau'} \intt{23} \ardl  \\
\tetav0 \ardl \intt{15} \ardr \tauv0 \\
         D_{0,\tau} \intt{35}  D_{\theta ,0} \\
\raisebox{-2mm}{\tauv0}
\ardr \intt{15} \ardl
\raisebox{-2mm}{\tetav0} \\
        D_{0,0}
\end{array}                      \label{e42a}
\end{equation}
        The points of the curve $ D_{(\theta ,0)} $ differ from that of
        $ D_{(0,0)} $ by a single commutator
        \begin{equation}
        [Y_{\pm}, X_{\pm}] = \theta (\mp e^{\pm tL} \pm e^{-hH}).
                                        \label{e42}
        \end{equation}
        For the contraction curve $ D_{(0,\tau)} $ the coproducts
        $ \Delta H, \Delta Y_{\pm} $ have the initial form
        (\ref{e34},\ref{e38})
        while the multiplication coincides with that of $D_{(0,0)}$
        (see (\ref{e39})).
        The vector field $ V(\theta) $ is described by the deformation
        function:
        \begin{equation}
        \begin{array}{rcl}
        \rule{1mm}{0mm} [X_+,X_-]_{(1)} & = & \theta^2 h^2
        \frac{e^{tL} - e^{-tL} }{1- e^{-th}} ,    \\
        \rule{1mm}{0mm} [H,X_{\pm}]_{(1)} & = & \pm \theta
        \frac{2th^2 }{1- e^{-th}} Y_{\mp} , \\
        \rule{1mm}{0mm} [Y_{\pm},X_{\mp}]_{(1)} & = &
        \pm \theta  h^2 Y^2_{\pm}
        +( e^{-th} - 1 ) X_{\mp} Y_{\pm},
        \end{array}                     \label{e43}
        \end{equation}
        \begin{equation}
        \begin{array}{rcl}
        (\Delta - \Delta^{\mbox{\scriptsize opp}})_{(1)} H  &
        = & \frac{2th^2}{1-e^{-th}}
        (Y_- \otimes Y_+ - Y_+ \otimes Y_-), \\
        (\Delta - \Delta^{\mbox{\scriptsize opp}})_{(1)} Y_{\pm} & = & \mp h^2
[
        (e^{-hH} \otimes Y_{\pm})(Y_{\mp} \otimes Y_{\pm}) -
        (Y_{\pm} \otimes e^{-hH})(Y_{\pm} \otimes Y_{\mp})].
        \end{array}                     \label{e44}
        \end{equation}
        In the vector field $ W(\tau) $ only some antisymmetric
        multiplication structure constants on $A \bigwedge A$
        are different from zero,
        \begin{equation}
        \begin{array}{rcl}
        \rule{1mm}{0mm} [H,X_{\pm}] & = &  \pm  \tau
        \frac{2th^2 }{1- e^{-th}} Y_{\mp} , \\
        \rule{1mm}{0mm} [Y_{\pm},X_{\mp}] & = &
        \pm \tau h^2 Y^2_{\pm},
        \end{array}                     \label{e45}
        \end{equation}
        \begin{equation}
        \rule{1mm}{0mm} [Y_{\pm},X_{\pm}]  =
        (\mp e^{\pm tL} \pm e^{-hH} ) .  \label{e46}
        \end{equation}
        Both $ V(\theta) $ and $ W(\tau) $ nontrivially depend on the
        co-ordinates of contraction curves. The limit values $ V(0) $
        and $ W(0) $  are defined by the relations (\ref{e44}) and
        (\ref{e46}) respectively.

        To complete the construction we must find such an
        $\varepsilon$-dependent family $ D_{\varepsilon} $
        of 2-dimen\-sion\-al varieties that in the limit $ \varepsilon
        \rightarrow 0 $ satisfies the canonical conditions (1)-(2).
        This may be done using the parameters $h$ and $t$.
        In the structure relations for $ D_{\theta ,\tau} $
        we perform the substitution:
        \begin{equation}
        h \Rightarrow \varepsilon^2, \, \, \,
        t \Rightarrow \varepsilon^2, \, \, \,
        X_{\pm} \Rightarrow \varepsilon X_{\pm}, \, \, \,
        Y_{\pm} \Rightarrow \varepsilon Y_{\pm}
                                        \label{e47}
        \end{equation}
        and obtain in the limit $ \varepsilon \rightarrow 0 $ the following
        compositions. For the internal points of $ D_{\varepsilon
         (\theta ,\tau)} $:
        \begin{equation}
        D_{ 0( \theta ,\tau)} \left\{ \begin{array}{rcl}
        \rule{1mm}{0mm} [X_+,X_-] & = & 2 \tau \theta^2 L, \\
        \rule{1mm}{0mm} [Y_{\pm},X_{\pm}] & = & - \theta (L \mp H), \\
        \Delta & -- & \mbox{primitive}.
                                        \end{array}
                                \right.
                                        \label{e48}
        \end{equation}
        and for the points of the boundary:
        \begin{equation}
        D_{0(\theta ,0)} \left\{ \begin{array}{rcl}
        \rule{1mm}{0mm} [Y_{\pm},X_{\pm}] & = & - \theta (L \mp H), \\
        \Delta & -- & \mbox{primitive},
                                        \end{array}
                                \right.
                                        \label{e49}
        \end{equation}
        \begin{equation}
        D_{0(0,\tau)} \left\{ \begin{array}{c}
        \mbox{Abelian}, \\
        \Delta  --  \mbox{primitive}.
                                        \end{array}
                                \right.
                                        \label{e50}
        \end{equation}
        \begin{equation}
        D_{0(0,0)} \left\{ \begin{array}{c}
        \mbox{Abelian}, \\
        \Delta  --  \mbox{primitive}.
                                        \end{array}
                                \right.
                                        \label{e51}
        \end{equation}

        All the Hopf algebras here are the classical universal
        enveloping  and the points of the curve $ D_{\varepsilon (0,\tau)} $
        are found to be trivialised. The deformation diagram obtains
        the form
\begin{equation}
\begin{array}{ccc}
& D_{\theta' ,\tau'}  &     \\
& \mid & \ardr \tauv0 \intt{3} \ardl \\
\mbox{\scriptsize $\theta \rightarrow 0$}
        & \mid & D_{\theta' ,0}    \\
&\downarrow & \ardl \tetav0 \intt{9} \\
& D_{0,0} & \approx  U(\mbox{Ab}) \intt{15}
        \end{array}                      \label{e52}
        \end{equation}
        It reflects the fact that the algebra $ D_{0(\theta ,\tau)} $
        -- a two-step classical first order deformation (of Abelian
        algebra) -- can also be treated as the second order deformation
        for each $ \tau \neq 0 $.

        So the limit $ \varepsilon \rightarrow 0 $ leads us to the
        degenerated case. The quantization removes this degeneracy.

        \section{Acknowledgments}
        I am heartily grateful to Prof. J.Lukierski for fruitful
        discussions and to all the scientists of the Institute of
        Theoretical Physics of Wroclaw University for their warm
        hospitality.

        The work is supported in part by the International Science
        Foundation, Grant N U9J000, and by the Russian Foundation
        for Fundamental Research, Grant N 95-01-00569a.

        \section{Appendix}
        The transformation $ B(t) $ and the reparametrization (\ref{e7})
        lead to the following 2-parametric variety
        $ Q(sl(n,{\bf C}),sl^*(n,{\bf C}))$ for quantum $sl(n,{\bf C})$
        algebras. (We use here the basis introduced in \cite{Ross}.)
        \begin{eqnarray*}
        \rule{1mm}{0mm}
        [H_i,H_j] &  = & 0, \rule{2cm}{0cm} i,j=1,\dots,n-1; \\
        \rule{1mm}{0mm}
        [H_i,X_{\pm(j,j+1)}] & = & \pm t \alpha_{ij} X_{\pm(j,j+1)}, \\
        \rule{1mm}{0mm}
        [X_{+(i,i+1)},X_{-(j,j+1)}]_{\exp(ht/2(\delta_{j+1,n} +
        \delta_{i+1,n}))} & = & t^2 \delta_{ij} e^{(ht/2)\delta_{i+1,n}}
        \frac{e^{-hH_{i;n-1}} - e^{-hH_{i+1;n-1}}}{e^{-ht/2} - e^{ht/2}}, \\
        \rule{1mm}{0mm}
        [X_{\pm (i,i+1)},X_{\pm (j,j+1)}]_{\exp(\pm ht/2(\delta_{i+1,n}
        - \delta_{j+1,n}))} & = & 0, \rule{2cm}{0cm} \mbox{for}
        \mid i-j \mid >1;
        \end{eqnarray*}
        \[
        \begin{array}{l}
        e^{\pm ht(\delta_{j,i+1} + \delta_{j+1,n})}(X_{\pm (i,i+1)})^2
        X_{\pm (j,j+1)} \\
        - e^{\pm ht/2(\delta_{j,i+1} + \delta_{j+1,i} +
        \delta_{j+1,n} + \delta_{i+1,n})}(e^{ht/2} - e^{-ht/2})
        X_{\pm (i,i+1)}X_{\pm (j,j+1)}X_{\pm (i,i+1)}  \\
        + e^{\pm ht(\delta_{j+1,i} + \delta_{i+1,n})}X_{\pm (j,j+1)}
        (X_{\pm (i,i+1)})^2 = 0, \rule{2cm}{0cm}
        \mbox{for} \mid i-j \mid =1;
        \end{array}
        \]
        \begin{eqnarray*}
        \Delta H_i & = & H_i \otimes {\bf 1} + {\bf 1} \otimes H_i, \\
        \Delta X_{\pm(i,i+1)} & = & X_{\pm(i,i+1)} \otimes
                e^{-(h/2)H_{i+1;n-1}} +
                e^{-(h/2)H_{i;n-1}} \otimes X_{\pm(i,i+1)}, \\
        S(H_i) & = & - H_i, \\
        S(X_{\pm (i,i+1)}) & = & -e^{hH_{i;n-1}/2}
                               X_{\pm (i,i+1)} e^{hH_{i+1;n-1}/2}, \\
        \varepsilon (X_{\pm (i,i+1)}) & = & \varepsilon (H_i) = 0
        \end{eqnarray*}
        Here $ \alpha_{i,j} $ is the Cartan matrix and
        \[
        H_{i,n-1} = H_{i,i+1} + H_{i+1,i+2} + \cdots + H_{n-2,n-1}.
        \]
        The modified $(\pm)$-co-ordinate functions \cite{Lya2}
        were applied so that the relation
        \[
        te^{\mp ht/2(1+ \delta_{i,j})}X_{\pm(i,j)} e^{-hH_{i+1,n-1}/2} =
        [X_{\pm (i,i+1)},X_{\pm (i+1,j)}]_{\exp(\pm (ht/2)\delta_{j,n})}.
        \]
        holds.
        In these terms the dual group $ SL^*(n,{\bf C}) $ has the simplest
        form. The coproduct for an arbitrary basic element $ X_{\pm(i,j)}$
        looks like
        \begin{eqnarray*}
        \Delta X_{\pm(i,j)} & = & X_{\pm (i,j)} \otimes e^{-hH_{j;n-1}/2}
                                + e^{-hH_{i;n-1}/2} \otimes X_{\pm (i,j)} \\
                        & + & \frac{(1-e^{\pm ht})}{t} \sum_{k=i+1}^{j-1}
                              X_{\pm(i,k)} \otimes X_{\pm(k,j)}.
        \end{eqnarray*}
        In the limit $ t \rightarrow 0 $ it describes the compositions
        of the solvable group $ SL^*(n,{\bf C}) $:
        \begin{eqnarray*}
        \Delta H_i & = & H_i \otimes {\bf 1} + {\bf 1} \otimes H_i, \\
        \Delta X_{\pm(i,j)} & = & X_{\pm (i,j)} \otimes e^{-hH_{j;n-1}/2}
                                + e^{-hH_{i;n-1}/2} \otimes X_{\pm (i,j)} \\
                        & \mp & h \sum_{k=i+1}^{j-1}
                              X_{\pm(i,k)} \otimes X_{\pm(k,j)}.
        \end{eqnarray*}

\end{document}